\documentclass[12pt,final,leqno]{amsart}
\RequirePackage{lmodern} 

\usepackage{amsfonts}
\usepackage{amsmath}
\usepackage{amssymb}
\usepackage{amsthm}
\usepackage{appendix}
\usepackage{changepage}
\usepackage{enumerate}
\usepackage{enumitem}
\usepackage{eucal}
\usepackage{latexsym}
\usepackage{mathrsfs}
\usepackage{scrextend}
\usepackage{textcomp}
\usepackage{todonotes}
\usepackage{verbatim}
\usepackage{pdflscape}
\usepackage{lscape}

\usepackage{xfrac}

\relax							

\newlength{\pardefault}							
\newenvironment{psmallmatrix}
  {\left(\begin{smallmatrix}}
  {\end{smallmatrix}\right)}

\newcommand{\innerProductReg}[2]{\left\langle #1 , #2 \right\rangle}
\newcommand{\innerProductTri}[3]{\left\langle #1 , #2 \right\rangle_{#3}}

\newcommand{\norm}[1]{\| #1  \|}
\newcommand{\normTwo}[2]{\| #1  \| _{#2}}

\newcommand{\real}{\mathfrak{Re} ~}
\newcommand{\imag}{\mathfrak{Im} ~}

\newcommand{\reall}[1]{\mathfrak{Re} \left( #1  \right)}
\newcommand{\imagg}[1]{\mathfrak{Im} \left( #1  \right)}

\newcommand\restr[2]{{
  \left.\kern-\nulldelimiterspace 
  #1 
  \vphantom{\big|} 
  \right|_{#2} 
  }}
\newcommand{\defEq}{\overset{\text{def} } {=}}

\newcommand{\tauBa}[1]{\tau{(#1)}}

\DeclareMathOperator*{\col}{Col}
\DeclareMathOperator*{\row}{Row}

\newcommand{\forceOddPage}{
\ifodd\theCurrentPage\pagebreak\fi}

\newcommand{\forceEvenPage}{
\ifodd\value{\theCurrentPage}\newpage\else\pagebreak\fi}

\newcommand*\cleartoleftpage{%
  \clearpage
  \ifodd\value{page}\hbox{}\newpage\fi
}

\newcommand*\cleartorightpage{%
  \clearpage
  \unless\ifodd\value{page}\hbox{}\newpage\fi
}

\newcommand{\ignore}[1]{}

\newtheorem{Pa}{Paper}[section]
\newtheorem{Tm}[Pa]{{\bf Theorem}}

\newtheorem{Cy}[Pa]{{\bf Corollary}}
\newtheorem{Rk}[Pa]{{\bf Remark}}

\newtheorem{theorem}[Pa]{{\bf Theorem}}
\newtheorem{lem}[Pa]{{\bf Lemma}}
\newtheorem{corollary}[Pa]{{\bf Corollary}}

\newtheorem{definition}[Pa]{{\bf Definition}}

\newenvironment{pf}[1][\unskip]{
\par
\noindent
\paragraph{{\bf Proof #1:}}
\noindent
}
{\hfill$\square$\\}

\def\R{\mathbb R}
\def\N{\mathbb N}
\def\C{\mathbb C}
\def\Z{\mathbb Z}
\def\T{\mathbb T}

\def\z{\zeta}
\def\e{\varepsilon}

\usepackage{hyperref}
\hypersetup{hypertexnames=false} 

\hypersetup{
  colorlinks,
  citecolor	= violet,
  linkcolor	= blue,
  urlcolor	= red}

% ...

\usepackage{tabularx}

\newcolumntype{P}[1]{>{\centering\arraybackslash}p{#1}}
\newcolumntype{M}[1]{>{\centering\arraybackslash}m{#1}}

\usepackage{chngcntr}
\counterwithin{equation}{section}

\begin{document}


\title
[Carath\'eodory functions in Riemann surfaces]
{Carath\'eodory functions on Riemann surfaces and reproducing kernel spaces}


\author[D. Alpay]{Daniel Alpay}
\address{(DA) Schmid College of Science and Technology,
Chapman University, One University Drive Orange, California 92866, USA}
\email{alpay@chapman.edu}

\author[A. Pinhas]{Ariel Pinhas}
\address{(AP) Department of mathematics,
Ben-Gurion University of the Negev, P.O. Box
653, Beer-Sheva 84105, Israel}
\email{arielp@post.bgu.ac.il}

\author[V. Vinnikov]{Victor Vinnikov}
\address{(VV) Department of mathematics,
Ben-Gurion University of the Negev, P.O. Box
653, Beer-Sheva 84105, Israel}
\email{vinnikov@math.bgu.ac.il}

\date{}

\thanks{The first author thanks the Foster G. and Mary McGaw Professorship in
Mathematical Sciences, which supported this research}


\begin{abstract}
Carath\'eodory functions, i.e. functions analytic in the open upper half-plane and with a positive real part there, play an important role 
in operator theory, $1D$ system theory and in the study of de Branges-Rovnyak spaces.
The Herglotz integral representation theorem associates to each Carath\'eodory function a 
positive measure on the real line and hence allows to further examine these subjects.
In this paper, we study these relations when the Riemann sphere 
is replaced by a real compact Riemann surface.
The generalization of Herglotz's theorem to the compact real Riemann surface setting is presented.
Furthermore, we study de Branges-Rovnyak spaces associated with functions with positive real-part defined 
on compact Riemann surfaces. Their elements are not anymore functions, but sections of a related line bundle.
\end{abstract}

\maketitle


\setcounter{tocdepth}{1}
\tableofcontents


\section{Introduction and overview}



A Carath\'eodory function $\varphi(z)$, that is, analytic with positive real part in the open upper half-plane ${\mathbb C}_+$,
admits an integral representation, also known by the Herglotz's representation theorem (see e.g \cite{MR48:904,RosenblumRovnyak}).
More precisely, the Carath\'eodory function $\varphi(z)$ can be written as:
\begin{equation}
\label{27-octobre-2000}
\varphi(z)=iA-iBz - i\int_{{\mathbb
R}}\left(\frac{1}{t-z}-\frac{t}{t^2+1} \right)d\mu(t),
\end{equation}
where $A\in{\mathbb R}$, $B\geq 0$ and $d\mu(t)$ 
is a positive measure on the real
line such that
$$\int_{\mathbb R}\frac{d\mu(t)}{t^2+1}<\infty.$$
One of the main two objectives of this paper is to
extend \eqref{27-octobre-2000} to the non zero genus case; this is done in Theorem \ref{caraTmRS}.
The second point is to extend \eqref{4-juin-2000} also to the non zero genus case, and this result is
presented in Theorem \ref{thm41}.
\smallskip

The right handside of \eqref{27-octobre-2000} defines an analytic extension
of $\varphi(z)$ to ${\mathbb C}\setminus{\mathbb R}$ such that
$\overline{\varphi(\overline{z})}+\varphi(z)=0$.
Thus, for any $z,w
\in{\mathbb C}\setminus {\mathbb R}$, we have
\begin{equation}
\label{4-juin-2000}
\frac{\varphi(z)+\overline{\varphi(w)}}{-i(z-\overline{w})}=
B+\int_{\mathbb R}\frac{d\mu(t)}{(t-z)(t-\overline{w})}=
B+\innerProductTri{\frac{1}{t-z}}{\frac{1}{t-w}}{{\bf L}^2(d\mu)},
\end{equation}
where ${\bf L}^2(d\mu)$ stands for the Lebesgue (Hilbert) space associated with the measure $d\mu$.
Thus, the kernel
$\frac{\varphi(z)+\overline{\varphi(w)}}{-i(z-\overline{w})}$ is
positive in ${\mathbb C}\setminus{\mathbb R}$. When $B=0$, the associated
reproducing kernel Hilbert space, denoted by $\mathcal{L}(\varphi)$, is described in the theorem below.
\begin{Tm}[{\cite[Section 5]{MR0229011}}]
\label{Thm21}
The space $\mathcal{L}(\varphi)$ consists of the functions of the form
\begin{equation}
\label{Thm21A}
F(z)=\int_{\mathbb R}\frac{f(t)d\mu(t)}{t-z}
\end{equation}
where $f\in{\bf L}^2(d\mu)$. Furthermore,
$\mathcal{L}(\varphi)$ is invariant under the resolvent-like
operators $R_\alpha$, where for $\alpha \in \C \setminus \R$, is given by:
\begin{equation}\label{liberte'}
(R_\alpha F)(z)=\frac{F(z)-F(\alpha)}{z-\alpha}.
\end{equation}

Finally, under the hypothesis $\int_{\mathbb R}d\mu(t)<\infty$, the elements of
$\mathcal{L}(\varphi)$ satisfy:
$$
\lim_{y\rightarrow\infty}F(iy)=0.$$
\end{Tm}

Moreover, the resolvent operator satisfies $R_\alpha = (M- \alpha I)^{-1}$,
where $M$ is the multiplication operator defined by
\[
M \, F(z) = z F(z) - \lim_{z \rightarrow \infty} z F(z).
\]
We note that $M$ corresponds, through \eqref{Thm21A}, to the operator of multiplication by $t$ in ${\bf L}^2(d\mu)$.
\smallskip

We provide here the outline of the proof in order to
motivate the analysis presented in the sequel in the compact real Riemann setting.

\begin{pf}[of Theorem \ref{Thm21}]
Let
$N\in{\mathbb N}$, $w_1,\ldots  , w_N \in{\mathbb C}\setminus{\mathbb R}$
and
$c_1 \cdots  c_N\in{\mathbb C}$.
Then,
\[
F(z)
\defEq
\sum_{j=1}^{N}
c_j \frac{\varphi(z)+\overline{\varphi(w_j)}}{-i(z-\overline{w_j})}
=
\int_{\mathbb R}\frac{d\mu (t)}{t-z}f(t)
\]
where
\[
f(t)=\sum_{j=1}^{N}\frac{c_j}{t-\overline{w_j}}\in{\bf L}^2(d\mu).
\]
In view of \eqref{4-juin-2000}, we have
\[
\|F\|^2_{\mathcal{L}(\varphi)}
=
\|f\|^2_{{\bf L}^2(d\mu)}
=
\sum_{\ell,j}\overline{c_\ell}
\frac{\varphi(w_\ell)+\overline{\varphi(w_j)}}{-i(w_\ell-\overline{w_j})}c_j.
\]
The first claim (Equation \ref{Thm21A}) follows by the fact that the linear span of the functions
$\frac{1}{-i(z-\overline{w})}$, $w\in{\mathbb C}\setminus{\mathbb R}$
is dense in ${\bf L}^2(d\mu)$.
\smallskip

Next, let $F(z)=\int_{\mathbb R}\frac{f(t)d\mu(t)}{t-z}\in\mathcal{L}(\varphi)$. Then,
\begin{equation}
\label{tatche}
(R_\alpha F)(z)=\int_{\mathbb R}\frac{f(t)d\mu(t)}{(t-\alpha)(t-z)}
\end{equation}
belongs to $\mathcal{L}(\varphi)$ since $f(t)/(t-\alpha)\in{\bf L}^2(
d\mu)$ where $\alpha$ is lying outside the real line.
\end{pf}

Furthermore, using \eqref{tatche}, the structure identity 

\begin{equation}
\label{strucId}
[R_\alpha f, g] - [f, R_\beta g] - (\alpha - \overline{\beta}) [R_\alpha f, R_\beta g] = 0,\quad \alpha,\beta\in\mathbb C\setminus\mathbb R,
\end{equation}
holds in the $\mathcal{L}(\varphi)$ spaces.  
In fact, this is an "if and only if" relation.
If the identity \eqref{strucId} holds in the space $\mathcal L$ of functions analytic in $\mathbb C\setminus\mathbb R$,
then $\mathcal{L}=\mathcal{L}(\varphi)$ for some Carath\'eodory function $\varphi(z)$ (see \cite[Theorem 6]{MR0229011}).
\smallskip

Using the observation that an ${\bf L}^2(d\mu)$ space is finite dimensional if and only if the measure $d\mu$ is has only a singular part, consisting of a finite number of jumps,
we may continue and mention the following result (see for instance \cite{dbbook}):

\begin{Tm}
\label{finiteDimentionalLphi}
Let $\varphi(z)$ be a Carath\'eodory function associated via \eqref{27-octobre-2000} to a positive measure $d \mu$
and let $\mathcal{L}(\varphi)$ be the corresponding 
reproducing kernel Hilbert space.
Then the following are equivalent:
\begin{enumerate}
\item $\mathcal{L}(\varphi)$ is finite dimensional.
\item ${\rm dim} \,	 {\bf L}^2(d \mu) < \infty$.
\item $d \mu$  is a jump measure with a finite number of jumps.
\item The Carath\'eodory function is of the form 
$$\varphi(z) = i A + i B z + \sum _{j=1}^N \frac{i c_j}{z-t_j},$$
where $c_j,B>0$, $A, t_j\in \R$ for all  $1\leq j\leq N$.
\end{enumerate}
\end{Tm}

\begin{Rk}
There are two different ways to obtain the positive measure $d\mu$ given in \eqref{27-octobre-2000}.
\begin{enumerate}
\item Using the Cauchy formula on the boundary and the Banach-Alaoglu Theorem.
\item Using the spectral theorem for $R_0$ in the space $\mathcal{L}(\varphi)$. 
In this case, $R_0$ is self-adjoint and the measure $d\mu$ is given 
by $d\mu(t) = \innerProductReg{dE(t)u}{u}$ where $E$ is the spectral measure of $R_0$.
\end{enumerate}
In this paper we focus on the first approach, while in \cite{AVP3}
we explore the second approach.
\end{Rk}

We mention here that  Carath\'eodory functions are the characteristic functions or transfer functions 
of selfadjoint vessel or impedance $2D$ systems, respectively.
Furthermore, they are also related to de Branges-Rovnyak spaces $\mathcal{L}(\varphi)$ 
of sections of certain vector bundles defined on compact Riemann surfaces of non zero genus.
These subjects and interconnections are further studied by the authors in \cite{AVP3}.\\

{\bf Outline of the paper:}
The paper consists of five sections besides the introduction.
In Section \ref{secPrel},
we give a brief overview of compact real Riemann surfaces and the associated Cauchy kernels. 
\smallskip

In Section \ref{secHerg}, we describe explicitly the Green function on $X$ in terms of the canonical homology basis. 
As a consequence, we present the Herglotz representation theorem for compact real Riemann surfaces. 
We utilize the integral representation of Carath\'eodory functions in order to study, in Section \ref{secdBLphi}, 
the de Branges space $\mathcal{L}(\varphi)$.
\smallskip

In Section \ref{secPhiSingleVal}, we examine the case where
$\varphi(z)$ is a single-valued function which defines a contractive function $s(z)$ through the Cayley transformation.
Hence, we may determine the relation between the de Branges spaces $\mathcal{L}(\varphi)$ and the the de Branges Rovnyak space $\mathcal H (s)$ associated to $s$.
Finally, in Section \ref{chSumm43}, we summarize some of the results by comparing the $\mathcal{L}(\varphi)$ theory 
in the genus zero case and the real compact Riemann surfaces of genus $g>0$.

\section{Preliminaries}
\label{secPrel}

In this section, we give a 
brief review of the basic properties and definitions of compact real Riemann surfaces.
We replace the open upper-half plane (or, more precisely,
its double, i.e the Riemann sphere) by a compact real Riemann surface $X$ of genus $g>0$.
\smallskip

A survey of the main tools required in the present study (including the prime form and the Jacobian) can be found in \cite[Section 2]{av3},
and in particular the descriptions of the Jacobian variety of a real curve and the real torii is in \cite{vinnikov5}.
For general background, we refer to \cite{fay1,GrHa,gunning2,mumford1} and \cite{mumford2}.
\smallskip

It is crucial to choose a canonical basis to the homology group $H_1(X,\mathbb Z)$ which is symmetric, 
in some sense, under the involution $\tau$
(for more details we refer to \cite{gross1981real}. Here we use the conventions as in \cite{av3,vinnikov5}).
Let $X_{\mathbb R}$ be the be set of the invariant points under $\tau$, $X_{\mathbb R} = \{ p\in X | \tauBa{p} = p\}$, 
which is always assumed to be not empty.
Then $X_{\mathbb R}$ consists of $k$ connected components denoted by $X_j$ where $j=0,...,k-1$
(disjoint analytic simple closed curves).
We choose for each component $X_j$ a point $p_j \in X_j$. 
Then, we set $A_{g+1-k+j} = X_j$ and $B_{g+1-k+j}=C_j - \tauBa{C_j}$, where $j=1,...,k-1$ and 
$C_j$ is a path from $p_0$ to $p_j$ which does not contain any other fixed point.
We can extend to the homology basis $A_1,...,A_g,B_1,...,B_g$ under which
the involution is given by $\begin{psmallmatrix}I & H\\0 & -I\end{psmallmatrix}$,
where the matrix $H$ is given by
\[
H=
\left( \begin{smallmatrix}
0 & 1 \\
1 & 0  \\
& & \ddots \\
& & & 0 & 1 \\
& & & 1 & 0  \\
& & & & & & 0  \\
& & & & & & & \ddots \\
& & & & & & & & 0 \\
\end{smallmatrix} \right)
\quad
{\rm and}
\quad
H=
\left( \begin{smallmatrix}
1 &  \\
  &  \ddots \\
  &   & 1 &  \\
  &   &   & 0 \\
  &   &   &   &  \ddots  \\
  &   &   &   &   & 0 \\
\end{smallmatrix} \right)
,
\]
for the dividing case and the non-dividing case, respectively.
In both cases, $H$ is of rank of $g+1-k$.
Then, we choose a normalized basis of holomorphic differentials on $X$ satisfying $\int _{A_i } \omega_j = \delta_{ij}$.
The matrix $Z \in \mathbb C ^{g \times g}$, with entries $Z_{i,j} = \int _{B_i } \omega_j $, is symmetric, with positive real part, satisfies 
\[
Z^* = H - Z
\]
and is referred as the period matrix of $X$ associated with the basis $\left( \omega_j \right)_{j=1} ^g$.
The Jacobian variety is defined by
$J(X) = \mathbb C ^ g \backslash \Gamma$, where $\Gamma = \mathbb Z ^g + Z \mathbb Z ^g$, 
and the Abel-Jacobi map from $X$ to the Jacobian variety is given by 
\[
\mu:p \rightarrow \begin{pmatrix}\int_{p_0}^p \omega_1\\ \vdots\\ \int_{p_0}^p \omega_g\end{pmatrix}.
\]
It is convenient to define 
\begin{equation}\label{eqZHY}Z= \frac{1}{2}H + i Y^{-1}.\end{equation}
\smallskip

We denote the universal covering of $X$ by $\pi:\widetilde{X}\rightarrow X$.
The group of deck transformations of $X$, denoted by the $\mathrm{Deck} (\widetilde{X} / X)$, consists of the homeomorphisms $\mathcal{T}: \widetilde{X} \rightarrow \widetilde{X}$
such that $\pi_X \circ \mathcal{T} = \pi_X$. 
It is well-known that the group of deck transformations on the universal covering is isomorphic to the fundamental group $\pi_1(X)$.
\smallskip

The analogue of the kernel $\frac{1}{-i(z-\overline{w})}$, is given by
$\frac{K_{\zeta}(u,\tauBa{v})}{-i}$ where
\[
K_{\zeta}(u,v) \defEq \frac{ \vartheta [{\zeta}](v-u)} { \vartheta [ \zeta ](0)E(v,u)}.
\]
The analogue of the kernel $\frac{1-s(z)\overline{s(w)}}{-i(z-\overline{w})}$
is now given by the expression
\begin{equation*}
K_{\tilde{\zeta},s}(u,v)=\frac{\vartheta [\tilde{\zeta}  ](\tauBa{v}-u)}
{i\vartheta [\tilde{\zeta}  ](0)E(u,\tauBa{v})}-
s(u)
\frac{ \vartheta [{\zeta}](\tauBa{v}-u)} {i\vartheta [ \zeta ](0)E(u,\tauBa{v})} 
\overline{s(v)}
,
\end{equation*}
where $u$ and $v$ are points on
$X$ (see \cite{vinnikov4} and \cite{vinnikov5}).
Furthermore,
$\zeta$ and $\tilde{\zeta}$
are points on the Jacobian $J(X)$ 
(in fact $\zeta$ and $\tilde{\zeta}$ belong to the real torii $T_\nu$, see \cite{MR1634421}) 
of $X$ such that $\vartheta(\zeta)$ and
$\vartheta(\tilde{\zeta})$ are nonzero and:

\begin{enumerate}

\item $\vartheta [\zeta ]$ denotes the theta function of $X$ with characteristic 
$\left[ \begin{array}{c} a \\ b \end{array} \right]$
where $\zeta=b+Za$ (with $a$ and $b$ in ${\mathbb R}^g$).
\item 
$E(u,v)$ is the prime form on $X$, for more details see \cite{fay1,mumford2}.
\item 
For fixed $v$, the map $u\mapsto K_{\widetilde{\zeta},s}(u,v)$ is a
multiplicative half order differential (with multipliers corresponding
to $\tilde{\zeta}$).
\item 
$s$ is a map of line bundles on $X$ with multipliers corresponding to $\tilde{\zeta}-\zeta$ and
satisfying $s(u)s(\tauBa{u})^*=1$.
\end{enumerate}

The analogue of the operators \eqref{liberte'} is given now by
\begin{equation*}
R_\alpha^{y}f(u)=
\frac{f(u)}{y(u)-\alpha}-\sum_{j=1}^n \frac{1}{ d y(u^{(j)})}
\frac{\vartheta[\zeta](u^{(j)}-u)}{\vartheta[\zeta](0)
E(u^{(j)},u)}
f(u^{(j)}),
\end{equation*}
where $y$ is a real meromorphic function of degree $n$ and $\alpha\in\C$ is
such that there are $n$ distinct points $u^{(j)}$ in $X$ such that $y(u^{(j)})=
\alpha$ and where $f$ is a section of
$L_{\zeta}\otimes\Delta$ analytic
at the points $u^{(j)}$. Furthermore, ({\cite[Lemma 4.3]{av3}})
the Cauchy kernels are eigenvectors of $R^y_\alpha$ with eigenvalues $\frac{1}{\overline{y(w)}-\alpha}$.

We conclude with the definition of the model operator, $M^{y}$ \cite[Equation 3-3]{MR1634421},
satisfying $(M^y - \alpha I ) ^{-1} = R_\alpha ^y$ for $\alpha$ large enough.
It is defined on sections of the line bundle $L_{\zeta}\otimes \Delta$ analytic 
at the neighborhood of the poles of $y$ and is explicitly given by
\begin{equation}
M^{y}f(u)
\label{m_y}
=
y(u)f(u) + \sum_{m=1}^{n}{c_m f(p^{(m)})
\frac
{\vartheta[\zeta](p^{(m)}-u)}
{\vartheta[\zeta] (0)E(p^{(m)},u)}},
\end{equation}
where $y(u)$ is a meromorphic function on $X$ with $n$ distinct simple poles, $p^{(1)},...,p^{(n)}$.


\section{Herglotz theorem for compact real Riemann surfaces}
\label{secHerg}

We first develop the analogue of Herglotz's formula for analytic functions with a positive real part in ${X}_+$,
instead of $\C_+$.
We consider the case of multi-valued functions but with purely imaginary period,
i.e. multi-valued functions that satisfy
\begin{equation*}
\varphi(\mathcal{T}(\widetilde{p}))=\varphi(\widetilde{p})+\chi(\mathcal{T}).
\end{equation*}
Here $\mathcal{T}$ is an element in the group of deck transformations 
on the universal covering of ${X}$ and
$$\chi:\,\,\pi_1(X)\rightarrow i{\mathbb R},$$
is a homomorphism of groups.
We call such a mapping an {\it additive function}.
Although in general it is not uniquely defined,
the real part of $\varphi(p)$ is well-defined.
\smallskip

The involution $\tau$ is extended on the universal covering of $X$.
In particular, for $\widetilde{x} \in \widetilde{X}$, an inverse under $\pi$ of an element in $X_\R$,
there exists $\mathcal{T}_{\widetilde{x}} \in \mathrm{Deck}(\widetilde{X} / X)$ such that 
$\tauBa{\widetilde{x}} = \mathcal{T} _{\widetilde{x}}(\widetilde{x})$.
Hence, since $\mathrm{Deck}(\widetilde{X} / X)$ is isomorphic to $\pi_1(X)$, we write $\mathcal{T}_{\widetilde{x}}$ in the form
$\mathcal{T}_{\widetilde{x}} = \sum_{j=1}^{g}{m_j A_j + n_j B_j}$, and we extensively use the notation 
\begin{align*}
n( \cdot )   : & \widetilde{X}  \longrightarrow \Z^g  \\
               & \widetilde{x} \longrightarrow (n_1 \cdots n_g)^{t}.
\end{align*}
We note that when $\widetilde{x} \in \widetilde{X}_0$ it follows that $n(\widetilde{x}) = 0$ and
$n(\widetilde{x}) = e_{g+1-k+j}$ whenever $\widetilde{x} \in \widetilde{X}_j$ for $j=1,...,k-1$
(where the set $e_1,...,e_g$ forms the canonical basis of $\R ^g$).
\begin{theorem}
\label{harmonicIntRep}
Let $X$ be a compact real Riemann surface and
let $\psi(p)$ be a positive harmonic function defined on ${X}\setminus {X}_{\mathbb R}$.
Then for every $p \in {X}\setminus {X}_{\mathbb R}$ there
exists a positive measure $d \eta(p,x)$ on $X_{\mathbb R}$ such that
\begin{align}
\label{la-guerre-commence1}
\psi(p)
 = &
\int_{X_{\mathbb R}}
\psi(x)
d \eta(p,x).
\end{align}
\end{theorem}

We start by presenting a preliminary lemma, revealing a useful property of the prime form.

\begin{lem}
\label{primeFormA}
Let $x$ be an element of $X_j$.
Then the prime form satisfies the following relation
\begin{align}
\label{primeFormEqA}
\overline{\frac{\partial}{\partial x}\ln E(\tauBa{p},x)}
& =
\frac{\partial}{\partial x}\ln E(p,x)-2\pi i \omega_{g-k-1+j}(x)
\\ & =
\frac{\partial}{\partial x}\ln E(p,x)-2\pi i \Big[ \omega_1(x)  \cdots  \omega_g(x) \Big] n(x).
\nonumber
\end{align}
\end{lem}
\begin{pf}
Let $x\in{ X}_j$, where $j=1,2,\ldots , k-1$.
Then, $\tau$ is lifted to the universal covering as follows
\[
\widetilde{x} - \tauBa{\widetilde{x}}
=
\sum_{\ell=1}^g m_\ell A_\ell + n_\ell B_\ell,
\]
where $n_\ell = \delta (g-k-1+j - \ell )$ 
and where $\delta$ stands for the Kronecker delta.
We also recall that the prime form (see \cite[Lemma 2.3]{av3}) satisfies
\begin{align}
\nonumber
E(\widetilde{p},\widetilde{u}_1)
= &
E(\widetilde{p},\widetilde{u}_2)
\exp
\left(
{-i \pi n^t \Gamma n + 2 \pi i (\widetilde{\mu}(\widetilde{p})-\widetilde{\mu}(\widetilde{u}_2))^t n}
\right) \times
\\&\times
\exp
\left(
2 \pi i (\beta^t_0 n - \alpha ^t _0 m )
\right)
\label{eqPrimeFormConj}
,
\end{align}
where $\widetilde{\mu}$ is the lifting of the Abel-Jacobi mapping to the universal covering
and where $\zeta = \alpha_0 + \beta_0 \Gamma$.
Thus, choosing $\tauBa{\widetilde{u}_2} = \widetilde{u}_1 = \widetilde{x}$, the relation in \eqref{eqPrimeFormConj} becomes
\begin{align}
\label{primeFormProp}
\ln \, E(\widetilde{p},\tauBa{\widetilde{x}})
= &
\ln \, E(\widetilde{p},\widetilde{x})
-i \pi \Gamma_{jj} + 2 \pi i (\widetilde{\mu}(\widetilde{p})-\widetilde{\mu}(\widetilde{x}))_j
+
\\ &+ 2 \pi i (\beta^t_0 n - \alpha ^t _0 m )
.
\nonumber
\end{align}
We note that by using \cite[Lemma 2.4]{av3}, the prime form satisfies the identity
$\overline{E(\tauBa{p},x)} = E(p,\tauBa{x})$.
It remains to differentiate \eqref{primeFormProp} with respect to $x$ and \eqref{primeFormEqA} follows.
\end{pf}

\begin{pf}[of Theorem \ref{harmonicIntRep}]
We show that the expression
\begin{align}
\label{la-guerre-commence2}
G(p,x) \defEq &
\pi \Big[ \omega_1(x)  \cdots  \omega_g(x) \Big] \cdot \left( \frac{n(x)}{2}   +  i  (Yp) \right)
-
\\ & \nonumber
-\frac{i}{2} \frac{\partial}{\partial x} \ln E(p,x)
,
\end{align}
is the differential with respect to $x$ of the Green function, where
$x\in X_\R$, $p \in X\setminus X_\R$ and where $\omega(x)$ is a section of the canonical bundle 
(denoted by $K_X$ and for an atlas $(V_j,z_j)$ defining the analytic structure of $X$, is given by cocycles $dz_j/dz_i$).
Here and in the following pages, with an abuse of notation, $Yp$ denotes $Y \widetilde{\mu}(\widetilde{p})$.
The existence of a Green function on a Riemann surface is a well-known result,
see for instance \cite[Chapter V]{bergman} or \cite[Chapter X]{Tsuji}.
Therefore, there exists a (unique) Green function, denoted by $g(p,x)$, with the
differential $G(p,x)$ which contains singularities of the form $\frac{1}{x-p}$
along its diagonal. Hence, it is enough to show that
the expression in \eqref{la-guerre-commence2} and the Green function 
satisfy the upcoming properties:
\begin{enumerate}
\item
{\it The function $g(x,p)$ contains a logarithmic singularity
while $G(x,p)$ has a simple pole at $p=x$.}
It follows immediately by using the prime form properties and moving to local coordinates 
that the following relation holds (see for instance \cite[Section II]{fay1}):
\[
\frac{i}{2} \frac{\partial}{\partial x} \ln E(p,x) 
= 
\frac{i}{2} \frac{\partial}{\partial v} \ln (t(u)-t(v)) 
= 
\frac{i}{2(t(u)-t(v))}
.
\]

\item
{\it The real part of the differential $G(x,p)$ is single-valued:}
Let $p$ and $p_1$ be two elements of $\widetilde{X}$ which are the pre-images of the same element in $X_j$, 
i.e. $\pi(p)=\pi(p_1) \in X_j $. It follows, using \eqref{eqZHY}, that 
\begin{equation}
\label{eqPp1}
\widetilde{\mu}(p) - \widetilde{\mu}(p_1) 
= 
n + \Gamma m
= 
n + \left(\frac{1}{2} H + i Y^{-1}\right) m,
\end{equation}
for some $n,m \in \R^g$ and thus, using again \cite[Lemma 2.3]{av3} we have, modulo $2\pi i$:
\begin{align*}
\ln\left( E(p_1,x) \right)
= &
\ln \bigg( E(p,x) \exp \big( 2 \pi i (\mu(x) - \mu(p))^t m \big) 
\times
\\ &  
\exp \left( 2 \pi i (\beta_0^t m - \alpha_0^t n)  -i \pi m^t \Gamma m  \right) \bigg)
\\ = &
\ln \, E(p,x) 
-
\frac{i}{2} \pi m^t H  m 
-  
\pi m^t  Y^{-1} m +
\\ & 
2 \pi i (\widetilde{\mu}(x) - \widetilde{\mu}(p))^t m + 2 \pi i (\beta_0^t m - \alpha_0^t n).
\end{align*}
Then, the real part of a multiplier of ${\it ln}\left( E(p,x) \right)$ is:
\begin{align}
\label{eqReLnMult}
\real
\big(
\ln E(p,x) -  &\ln E(p_1,x)
\big)
= \\ & \nonumber
2 \pi \left(
\frac{1}{2} m^t Y^{-1}  + \imagg {  \mu(p) - \mu(x)} ^t
\right) m.
\end{align}
Clearly, using \eqref{eqPp1}, we have that
$$ m = Y \, \imagg { \widetilde{\mu}(p) - \widetilde{\mu}(p_1)}$$ 
and hence
the derivative with respect to $x$ of \eqref{eqReLnMult}, is
\begin{align*}
\frac{\partial}{\partial x}
\real 
(\ln \, \left( E(p,x)\right) & - \ln \, \left( E(p_1,x) \right)
)
\\ = &
- 2 \pi \, \imag
\Big[\omega_1(x)  \cdots  \omega_g(x) \Big]
m
\\ = &
- 2 \pi \, \imag
\Big[\omega_1(x)  \cdots  \omega_g(x) \Big]
Y \widetilde{\mu} (p - p_1)
.
\end{align*}
Hence, $G(x,p)$ has the appropriate singularity and has a single-valued real-part if it is of the following form:
\[
\frac{\partial}{\partial x}
\ln\left( E(p,x) \right)+ 2 \pi \Big[\omega_1(x)  \cdots  \omega_g(x) \Big] 
Y \widetilde{\mu} (p) + h(x),
\]
for some $h(x)$ with purely imaginary periods.

\item
{\it The real part of the complex Green function vanishes on the boundary components:}
Let $x \in X_j$ for some $0 \leq j \leq k-1$ and let $p \in X_l$ for some $0 \leq l \leq k-1$ such that $p \neq x$.
Then, we integrate $G(x,p)$ with respect to $x$
and note that the integration of the vector $\Big[ \omega_1(x)  \cdots  \omega_g(x) \Big]$ is
just the Abel-Jacobi mapping at $x$.
Then, the Green function is:
\[
g(x,p)= \left( \frac{n(x)}{2}   +  i  (Yp) \right) \mu(x) -
\frac{i}{2} \ln E(p,x).
\]
We use the equality, see \cite[Lemma 2.4]{av3},
$$ E(x,p) = \overline{E(\tauBa{x},\tauBa{p})} $$
to conclude that whenever $x$ and $p$ are both real, the relation
\begin{align*}
\overline{ \frac{\partial}{\partial x} \ln{(E(x,p))}} = &
\frac{\partial}{\partial x} \ln{(E(\tauBa{x},\tauBa{p}))} \\ = &
\frac{\partial}{\partial x} \ln{(E(x,p))} + 2\pi i \omega _{g-k+j-1}(x)
\end{align*}
holds. Hence, the real part of $g(x,p)$, using Equation \ref{primeFormEqA}, is equal to
\begin{align*}
\reall{g(x,p)}
= &
i  (Yp) \mu(x) -
\frac{i}{2} \ln E(p,x)
+
\overline{i  (Yp) \mu(x) }
-
\\ &
\overline{\frac{i}{2} \ln E(p,x)}
+
\reall{h(x)}
\\ = &
\real
\frac{i}{2}
\left(
\ln E(\tauBa{p},\tauBa{x})
-
\ln E(p,x)
\right)
+
\reall{h(x)}
\\ = &
\frac{\pi }{2}
\omega(x) n(x)
+
\reall{h(x)}
,
\end{align*}
and therefore, setting $\real h(x) = - \frac{\pi}{2} w(x) n(x)$, the Green function vanishes on the real points.
\end{enumerate}
Thus, $G(x,p)$ is the differential of the complex Green function
and so, for any $p \in X \setminus X_\R$ and for
sufficient small $\e$, defines the solution to the Dirichlet problem, i.e.
\[
\psi(p)  = \int_{X_\R(\e)} \psi(x) G(x,p).
\]
Here, the integration contour is a collection of smooth simple closed curves 
located within a distance $\e$ approximating $X_\R$. 
We then consider a sequence $(\e_n)_{n \in \N}$ such that $\e_n \rightarrow 0$ as $n \rightarrow \infty$.
Then, by the Banach-Alaoglu Theorem (see for instance, \cite[p. 223]{MR1681462}),
there exists a subsequence $(\e_{n_k})_{k \in \N}$ such that the limit
$$\lim_{k \rightarrow \infty} \int_{X_\R(\e_{n_k})} \psi(x) G(x,p)$$
exists.
Thus, the weak-star limit defines a positive measure on $X_\R$ satisfying \eqref{la-guerre-commence1}.
\end{pf}

Using the previous result, we may state the Herglotz theorem for
real compact Riemann surfaces.

\begin{Tm}
\label{caraTmRS}
Let $X$ be a compact real Riemann surface of dividing-type.
Then an additive function $\varphi(x)$ analytic in
${X}\setminus {X}_{\mathbb R}$ with positive real part in ${X}\setminus {X}_{\mathbb R}$ and,
furthermore, satisfies
\begin{equation*}
\varphi(p)+\overline{\varphi(\tauBa{p})}=0,\quad p\in{ X}\setminus { X}_{\mathbb R},
\end{equation*}
if and only if
\begin{align}
\nonumber
\varphi(p)
 = &
\frac{\pi}{2}
\int_{X_{\mathbb R}}
\Big[\omega_1(x)  \cdots  \omega_g(x) \Big]
n(\widetilde{x}) \, \frac{d \eta(x)}{\omega(x)}
-
\frac{i}{2}\int_{X_{\mathbb R}}  \frac{\partial}{\partial x} \ln E(p,x) \, \frac{d \eta(x)}{\omega(x)} +
\\
\label{la-guerre-commence}
&
\pi i\int_{X_{\mathbb R}} \Big[ \omega_1(x)   \cdots  \omega_g(x) \Big] (Yp) \, \frac{d \eta(x)}{\omega(x)}
+iM.
\end{align}
Here,
$M$ is a real number,
$d \eta$ is a positive finite measure on $X _{\mathbb R}$,
$\omega(x)$ is a section of the canonical line bundle 
which is positive with respect to the measure $d \eta$.
\end{Tm}
\begin{pf}
We start with the "if" part as we compute $\overline{\varphi(\tauBa{p})}$:
\begin{eqnarray*}
\overline{\varphi(\tauBa{p})}
&=&
\frac{\pi}{2}
\int_{{X}_{\mathbb R}}
[\omega_1(x)   \cdots  \omega_g(x) ] n(\widetilde{x})
\frac{d \eta(x)}{\omega(x)}
-
\\
& &
\pi i
\int_{{X}_{\mathbb R}}
[\omega_1(x)   \cdots  \omega_g(x)] (\overline{Y p})~
\frac{d \eta(x)}{\omega(x)}
+
\\
& &
\frac{i}{2}
\int_{{X}_{\mathbb R}} \frac{\partial}{\partial x} \ln \overline{E(\tauBa{p},x)} \frac{d \eta(x)}{\omega(x)}
-iM
.
\end{eqnarray*}
Thus, using Lemma \ref{primeFormA} and since $\omega$ is real (i.e. $\overline{\tauBa{ \omega_i}} = \omega_i$),
we have:
\begin{align}
\nonumber
\overline{\varphi(\tauBa{p})}
=&
\frac{\pi}{2}
\int_{{X}_{\mathbb R}}
[\omega_1(x) \, \cdots \, \omega_g(x)] n(\widetilde{x})
\frac{d\eta(x)}{\omega(x)} -
\\
\nonumber
&
\pi i \int_{{ X}_{\mathbb R}}
[\omega_1(x) \, \cdots \, \omega_g(x)]
(Yp)~
\frac{d \eta(x)}{\omega(x)}
+
\\
\label{la-guerre-commence11}
&
\frac{i}{2}
\int_{{X}_{\mathbb R}} \frac{\partial}{\partial x} \ln E(p,x)
\frac{d \eta(x)}{\omega(x)}
-
\pi \sum_{j=1}^{k-1}\int_{{X}_j}\omega_j(x)
\frac{d \eta(x)}{{\omega(x)}}
-iM.
\end{align}
Summing up \eqref{la-guerre-commence} and \eqref{la-guerre-commence11}, leads to
\begin{align*}
\varphi(p)+\overline{\varphi(\tauBa{p})}
= &
\pi
\int_{{X}_{\mathbb R}}
[ \omega_1(x) \, \cdots \, \omega_g(x)] n (\widetilde{x})
\frac{d \eta(x)}{\omega(x)} -
\\
&
\pi \sum_{j=1}^{k-1}\int_{{X}_j} \omega_j(x)
\frac{d \eta(x)}{\omega(x)}
= 0.
\end{align*}
For the "only if" statement:
The real part of $\varphi(p)$
is positive, harmonic and with a single-valued real part
in $X \setminus X_{\mathbb R}$.
Thus, by Theorem \ref{harmonicIntRep}, $\real{ \varphi(p)}$
has an integral representation as given in \eqref{la-guerre-commence1}
for some positive measure $d \, \eta_{\varphi}$ on $X_\R$.
Finally, it is well-known that two analytic functions defined 
on a connected domain with the same real part
differ only by some imaginary constant.
Hence we may summarize that
\[
\varphi(p)
=
\int_{X_{\mathbb R}}
G(p,x) d \nu_{\varphi}(x)
+
iM,
\]
for some $M \in \R$.
\end{pf}

In the case where $X = \mathbb P^1$ coupled with the anti-holomorphic involution $z \rightarrow \overline{z}$,
we set $\omega = \frac{d \, t}{t^2 + 1}$ and then \eqref{27-octobre-2000} can be extracted from \eqref{la-guerre-commence} by setting,
\begin{align*}
d \nu (t) & = \frac{1}{2} d \eta (t) (t^2 +1), \qquad
B = \frac{1}{2} \eta (\infty),
\\
A & = M - \frac{1}{2} \int_{I} t \, d \eta (t) + \frac{1}{2} \int_{\mathbb R \backslash I} \frac{d \eta (t)}{t},
\end{align*}
where $I$ is any interval of $\mathbb R$ containing zero.
\smallskip

Similarly, in the case of the torus,
one may deduce H. Villat's formula, see \cite{MR1629812}.
(Akhiezer in \cite[Section 56]{MR1054205} presented a different but equivalent formula). 

\section{de Branges \texorpdfstring{$\mathcal{L}(\varphi)$}{ $\mathcal{L}(\varphi)$ } spaces in the nonzero genus case}
\label{secdBLphi}

In this section, we further study the reproducing kernel Hilbert space associated with
an additive function defined on a real compact Riemann space.
To do so, we utilize the Herglotz's integral representation proved inthe previous section in order to examine
$\mathcal{L}(\varphi)$ spaces and their properties.
First, we introduce the analogue of formula \eqref{4-juin-2000}.
\begin{Tm}
\label{thm41}
Let $X$ be a compact real Riemann surface of dividing type, $\z \in T_0$ and
let $\varphi$ be an analytic with positive real part in $X_+$.
Then, the identity
\begin{align}
\nonumber
\int_{{X}_{\mathbb R}}
&
\frac{\vartheta[\zeta](p-x)}{\vartheta[\zeta](0)E(x,p)}
\frac{\vartheta[\zeta](x-{\tauBa{q}})}{\vartheta[\zeta](0)E(x,{\tauBa{q}})}
\frac{d \eta(x)}{\omega(x)}
=
\frac{\vartheta[\zeta](p-{\tauBa{q}})}{\vartheta[\zeta](0)E({\tauBa{q}},p)}
\times
\\
\nonumber
&
\times\bigg[
\left({\varphi}(p)+\overline{{\varphi}(\tauBa{q})}\right)
+
\sum_{j=0}^{k-1}
a_{jj}\frac{\partial}{\partial z_j}
\ln\frac{\vartheta(\zeta)}{\vartheta(\zeta+p-{\tauBa{q}})}
-
\\
&
-2\pi i
\row_{i=0,\ldots,g} 
\left( 
\sum_{j=0}^{k-1}
a_{ji}
\right)
Y(p-{\tauBa{q}})
,
\label{jfk-le-26-octobre-2000}
\end{align}
\label{la-fin-du-sionisme?}
holds where
\begin{equation}
\label{a_j}
a_{ji} \defEq \int_{{X}_{j}}
\frac{\omega_i(x)}{\omega(x)}d \eta(x),\quad j=0,\ldots,k-1, \quad i=1,\ldots,g.
\end{equation}
\end{Tm}
Before heading to prove Theorem \ref{thm41}, we make a number of remarks.
The left hand side of \eqref{jfk-le-26-octobre-2000} may be written as
$$
\innerProductTri{K_{\zeta}(p,x)}{K_{\zeta}(x,\tauBa{q})}
{{\bf L}^2\left( X_{\mathbb R} , L_{\zeta} \otimes \Delta , \frac{d \eta(x)}{w(x)} \right) }
$$
and hence, it is precise the counterpart of the right hand side of \eqref{4-juin-2000}.
Furthermore, whenever we additionally assume zero cycles along the boundary components, that is,
\begin{equation}
\label{eqCycles}
\int_{{ X}_j}\frac{\omega_i(x)}{\omega(x)}d \eta(x)=0, \quad j=0,\ldots,k-1,
\end{equation}
the right hand side of \eqref{jfk-le-26-octobre-2000}
is the counterpart of the left hand side of \eqref{4-juin-2000}.
Hence, we may summarize and present the following result.
\begin{corollary}
\label{kernelAj0}
Let $X$ be a compact real Riemann surface of dividing type, let $\z \in T_0$ and
let $\varphi$ be an additive function on $X$ 
such that \eqref{eqCycles} holds. Then, the identity
\begin{align}
\left({\varphi}(p)+\overline{{\varphi}(\tauBa{q})}\right)
&
\frac{\vartheta[\zeta](p-{\tauBa{q}})}{\vartheta[\zeta](0)E({\tauBa{q}},p)}
\nonumber = \\ &
\innerProductTri{K_{\zeta}(p,u)}{K_{\zeta}(u,\tauBa{q})}
{{\bf L}^2\left( X_{\mathbb R} , L_{\zeta} \otimes \Delta , \frac{d \eta(x)}{w(x)} \right) }
,
\label{jfk-le-26-octobre-2000_SV}
\end{align}
holds.
\end{corollary}
\begin{pf}[of Theorem \ref{la-fin-du-sionisme?}]
Using the Herglotz-type formula \eqref{la-guerre-commence}, we may write
\begin{align}
\nonumber
(\varphi(p) + &\overline{\varphi(\tauBa{q})}) 
\frac{\vartheta[\zeta](p-{\tauBa{q}})}{\vartheta[\zeta](0)E({\tauBa{q}},p)}
= \\
=&
2\pi i\left(\int_{X_{\mathbb R}}
[ \omega_1(x)  \cdots  \omega_g(x) ] 
Y(p-\tauBa{q}) \frac{d\eta(x)}{\omega(x)}\right)
\frac{\vartheta[\zeta](p-{\tauBa{q}})}{\vartheta[\zeta](0)E({\tauBa{q}},p)}
\nonumber
- \\
\label{Thm43A}
&
\frac{i}{2}
\int_{X_{\mathbb R}}
\left(\frac{\partial}{\partial x}\ln E(p,x)-\frac{\partial}{\partial x}\ln E({\tauBa{q}},x)\right)
\frac{d\eta(x)}{\omega(x)}
\frac{\vartheta[\zeta](p-{\tauBa{q}})}{\vartheta[\zeta](0)E({\tauBa{q}},p)}.
\end{align}
Furthermore, by \cite[Proposition 2.10, p. 25]{fay1}, we have
\begin{align}
\frac{\partial}{\partial x}\ln\frac{E(p,x)}{E({\tauBa{q}},x)}
+ &
\sum_{j=1}^g\left(\frac{\partial}{\partial z_j}\ln \vartheta(\zeta+p-{\tauBa{q}})
-\frac{\partial}{\partial z_j}\ln \vartheta(\zeta)\right)\omega_j(x)
=
\nonumber
\\
\label{eq123}
&
\frac{E({\tauBa{q}}, p)}{E(x,{\tauBa{q}})E(x,p)}
\frac
{\vartheta[\zeta](x-{\tauBa{q}})\vartheta[\zeta](p-x)}
{\vartheta[\zeta](p-{\tauBa{q}})\vartheta[\zeta](0)}.
\end{align}
Thus, multiplying both sides of \eqref{eq123} by
$\frac{\vartheta[\zeta](p-{\tauBa{q}})}{\vartheta[\zeta](0)E({\tauBa{q}},p)}$,
leads to
\begin{align}
\nonumber
&
\frac{\vartheta[\zeta](p-{\tauBa{q}})}{\vartheta[\zeta](0)E({\tauBa{q}},p)}
\frac{\partial}{\partial x}\ln\frac{E(p,x)}{E({\tauBa{q}},x)}
=
\frac{\vartheta[\zeta](x-{\tauBa{q}})}{
\vartheta[\zeta](0)E(x,{\tauBa{q}})}
\frac{\vartheta[\zeta](p-x)}{
\vartheta[\zeta](0)E(x,p)}-
\\
&
\sum_{j=1}^g
\frac{\partial}{\partial z_j}
\left(
\ln \vartheta(\zeta+p-{\tauBa{q}})
-
\ln \vartheta(\zeta)\right)\omega_j(x)
\frac{\vartheta[\zeta](p-{\tauBa{q}})}{
\vartheta[\zeta](0)E({\tauBa{q}},p)}.
\label{Thm43B}
\end{align}
Finally, by substituting \eqref{Thm43B} into \eqref{Thm43A}, we conclude that the identity
\begin{align*}
(\varphi(p)-\overline{\varphi(\tauBa{q})}) &
\frac{\vartheta[\zeta](p-{\tauBa{q}})}{\vartheta[\zeta](0)E({\tauBa{q}},p)}
=
\frac{\vartheta[\zeta](p-{\tauBa{q}})}{
\vartheta[\zeta](0)E({\tauBa{q}},p)}
\times
\\
&
\bigg[
\frac{i}{2}
\sum_{j=1}^{g}
\int_{X_{j}}
\frac{\omega_j(x) d \eta(x)}{\omega(x)}
\frac{\partial}{\partial z_j}
\left(
\ln \vartheta(\zeta+p-{\tauBa{q}})
-\ln \vartheta(\zeta)
\right) +
\\
&
2\pi i\int_{X_{\mathbb R}}
\frac{[\omega_1(x)\,\cdots  \, \omega_g(x)]}{\omega(x)}Y(p-{\tauBa{q}})d \eta(x)
\bigg]
-
\\ &
\frac{i}{2}
\int_{{X}_{\mathbb R}}
\frac{\vartheta[\zeta](p-x)}{\vartheta[\zeta](0)E(x,p)}
\frac{\vartheta[\zeta](x-{\tauBa{q}})}{
\vartheta[\zeta](0)E(x,{\tauBa{q}})}\frac{d \eta(x)}{\omega(x)}
\end{align*}
follows.
Setting $a_j$ as in \eqref{a_j}, completes the proof.
\end{pf}

From this point and onward we assume that \eqref{eqCycles} holds.

\begin{definition}
Let $\varphi(x)$ be analytic in
${X}\setminus {X}_{\mathbb R}$ with positive real part in ${X}\setminus {X}_{\mathbb R}$.
The reproducing kernel Hilbert space of sections of the line bundle 
$L_\zeta \otimes \Delta$ with the reproducing kernel
\[
K(p,q)
=
(\varphi(p) + \overline{\varphi(\tauBa{q})})
\frac{\vartheta[\zeta](p-{\tauBa{q}})}{\vartheta[\zeta](0)E({\tauBa{q}},p)},
\]
is denoted by $\mathcal{L}(\varphi)$.
\end{definition}
The analogue of the first part of Theorem \ref{Thm21} is given below in Theorem \ref{phiIntPresentation}.
However, we first present a preliminary lemma that is required during this section
(see \cite[Ex. 6.3.2]{capb2}, in the unit-disk case).
\begin{lem}
\label{denseL2}
Let $X$ be a compact real Riemann surface of dividing type.
Then the linear span of Cauchy kernels
{\allowbreak
$\frac{\vartheta[\zeta](x-u)}{i \vartheta[\zeta](0)E(u,x)}$
}
where $u$ varies in $ X \setminus X_{\mathbb R}$ is dense in
{\allowbreak
${\bf L}^2\left( X_{\mathbb R} , L_{\zeta} \otimes \Delta , \frac{d \eta(x)}{w(x)} \right)$}.
\end{lem}
\begin{pf}
Let us assume that a section $f$ of $L_{\zeta} \otimes \Delta$ satisfies
\begin{equation}
\label{eqCkDense}
\int_{{X}_{\mathbb R}}K_{\zeta}(u, x)f(x)\frac{d \eta(x)}{\omega(x)} = 0,
\end{equation}
for all $u\in X \setminus X_\R$.
We recall that by \cite{av2}, there exists an isometric isomorphism from
${\bf L}^2\left( X_{\mathbb R} , L_{\zeta} \otimes \Delta , \frac{d \eta(x)}{w(x)} \right)$
to
$ ({\bf L}^2(\T))^n$ and therefore there exists an orthogonal decomposition, see \cite[Equation 4.14]{av2},
\begin{align*}
{\bf L}^2\bigg( & X_{\mathbb R} , L_{\zeta} \otimes \Delta , \frac{d \eta(x)}{w(x)} \bigg)
\\
= &
{\bf H}^2\left( X_{+} , L_{\zeta} \otimes \Delta , \frac{d \eta(x)}{w(x)} \right)
\oplus
{\bf H}^2\left( X_{-} , L_{\zeta} \otimes \Delta , \frac{d \eta(x)}{w(x)} \right)
.
\end{align*}
Furthermore, Equation \ref{eqCkDense}, for $u \in X_+$ is just the projection from
${\bf L}^2 \bigg( X_{\mathbb R} , L_{\zeta} \otimes \Delta , \frac{d \eta(x)}{w(x)} \bigg)$
into
${\bf H}^2 \left( X_{+} , L_{\zeta} \otimes \Delta , \frac{d \eta(x)}{w(x)} \right)$.
Thus, $P_+(f)(u)=0$ and, similarly, $P_-(f)(u)=0$
and we may conclude that $f=0$ and the claim follows.
\end{pf}
\begin{Tm}
\label{phiIntPresentation}
The elements of $\mathcal{L}(\varphi)$ are of the form
\begin{equation}
\label{l-phi}
F(u) = \int_{{X}_{\mathbb R}}K_{\zeta}(u, x)f(x)\frac{d \eta(x)}{\omega(x)},
\end{equation}
where $f(x)$ is a section of $L_\zeta\otimes\Delta$ which is square summable with respect to $\frac{d \eta (x)}{\omega(x)}$.
\end{Tm}
\begin{pf}
Equation \ref{l-phi} follows by Corollary \ref{kernelAj0}.
Let us set $N\in{\mathbb N}$,
then for any choice of $w_1,\ldots,w_N \in{X}\setminus X_{\mathbb R}$
and $c_1 \cdots  c_N \in {\mathbb C}$, the identity
\begin{align}
\label{lPhinorm}
F(u)
\defEq
&
\sum_{j=1}^{n}
c_j (\varphi(u) + \overline{\varphi(w_j)}) K_{\zeta}( u, w_j)
\\
=
&
\int_{X_{\mathbb R}} K_{\zeta}(u, x) f(x) \frac{d \eta (x)}{w(x)}
\nonumber
\end{align}
holds, where
\[
f(u)
=
\sum_{j=1}^{n} c_j K_{\zeta}( w_j, u)
\in
{{\bf L}^2\left( X_{\mathbb R} , L_{\zeta} \otimes \Delta , \frac{d \eta(x)}{w(x)} \right)}
.
\]
Due to Lemma \ref{denseL2}, the linear span of the kernels \eqref{jfk-le-26-octobre-2000_SV} is dense in
${\bf L}^2\left( X_{\mathbb R} , L_{\zeta} \otimes \Delta , \frac{d \eta(x)}{w(x)} \right)$
and hence \eqref{l-phi} follows.
\end{pf}
\begin{Tm}
The norm of an element $F$ in $\mathcal{L}(\varphi)$ is given by
\begin{equation*}
\normTwo{ F }{\mathcal{L}(\varphi)}
\defEq
\normTwo{f}{{\bf L}^2 \left( X_{\mathbb R} , L_{\zeta} \otimes \Delta , \frac{d \eta(t)}{w(t)} \right)}.
\end{equation*}
\end{Tm}
\begin{pf}
Since, by Lemma \ref{denseL2}, the linear span of the kernels \eqref{jfk-le-26-octobre-2000_SV} is dense in ${\bf L}^2(d \eta)$,
it is enough to check the eqaulity of the norms for a linear combination of the Cauchy kernels.
The norm of an element in the reproducing kernel Hilbert space $\mathcal{L}(\varphi)$ is given by
\[
\norm{F}^2_{\mathcal{L}(\varphi)}
=
\sum_{\ell,j=1}^{n}
\overline{c_\ell} (\varphi(u_\ell) + \overline{\varphi(u_j)}) K_{\zeta}( u_\ell, u_j)c_j.
\]
Then, by \eqref{lPhinorm}
\[
\norm{F}^2_{\mathcal{L}(\varphi)}
=
\sum_{\ell,j=1}^{n}
\overline{c_j}
\innerProductTri
{K_{\zeta}( u_\ell, w )}
{K_{\zeta}( w, u_j)}
{{\bf L}^2\left( X_{\mathbb R} , L_{\zeta} \otimes \Delta , \frac{d \eta(t)}{w(t)} \right)}
c_\ell,
\]
which is exactly the norm of $f(w)$ in
${\bf L}^2\left( X_{\mathbb R} , L_{\zeta} \otimes \Delta , \frac{d \eta(t)}{w(t)} \right)$.
\end{pf}

As an immediate consequence, whenever $y$ is a real function,
we may state an additional result.
It follows that $M^y$ is simply the multiplication operator in ${\bf L}^2(d \eta)$.

\begin{Tm}
Let $y$ be a meromorphic function with simple poles such that the poles of $y$ do not belong to the support of the measure $d \eta$.
Then, the multiplication model operator $M^{y}$, defined on $\mathcal{L}(\varphi)$, satisfies the following properties:
\begin{enumerate}
\item $M^{y}$ is given explicitly by
\begin{equation}
\label{MyLPhi}
\left(M^y F \right) (u)
=
\int_{{X}_{\mathbb R}} 
{ K_{\zeta}( u, x
)}f(x)y(x)\frac{d \eta(x)}{\omega(x)},
\end{equation}
where $f$ is a section of $L_\zeta\otimes\Delta$, which is square summable with respect to $d \eta$.
\item $\mathcal{L}(\varphi)$ is invariant under $M^{y}$.
\item $M^{y}$ is bounded.
\end{enumerate}
\end{Tm}
\begin{pf}
Considering the model operator \eqref{m_y} together with Theorem \ref{phiIntPresentation}, we conclude the following:
\begin{align*}
(M^y F)(u) = & y(u)F(u) + \sum_{m=1}^{n}{c_m F(p_m) K_{\zeta}( u, p_m)}
\\
= &
y(u)F(u) +  \sum_{m=1}^{n}{c_m \int_{X_{\mathbb R}} f(x) K_{\zeta} (p_m,x) \frac{d \eta(x)}{\omega(x)}  K_{\zeta} ( u, p_m)}
\\
= &
y(u)F(u) + \int_{X_{\mathbb R}} f(x)  \frac{d \eta(x)}{\omega(x)}  \sum_{m=1}^{n}{c_m K_{\zeta}( p_m,x) K_{\zeta} ( u, p_m)},
\end{align*}
where $p^{(1)},...,p^{(n)}$ are the distinct poles of $y$.
Using the collection formula \cite[Proposition 3.1]{av2} and using again Theorem \ref{phiIntPresentation}, we have:
\begin{align*}
(M^y F)(u)
= &
y(u)F(u) + \int_{X_{\mathbb R}} f(x)  \frac{d \eta(x)}{\omega(x)}  K_{\zeta}  (u,x) (y(x)-y(u))
\\
= &
\int_{{X}_{\mathbb R}} 
{ K_{\zeta}( u, x
)}f(x)y(x)\frac{d \eta(x)}{\omega(x)}.
\end{align*}
We note that $f$ is a section of $L_\zeta\otimes\Delta$ and remains so after multiplication by a meromorphic function $y$.
Furthermore, it is square summable with respect to the measure $d \eta$, since,
by assumption, the poles of $y$ lie outside the support of $d \eta$.
\end{pf}

\begin{corollary}
Let $y$ be a real meromorphic function on $X$ such that the poles of $y$
do not belong to the support of the measure $d \eta$.
Then, the multiplication model operator $M^{y}$ is selfadjoint.
\end{corollary}
\begin{pf}
The model operator $M^y$ satisfies $\innerProductReg{M^y F}{G} = \innerProductReg{ F}{M^yG}$ as follows from
\[
\innerProductReg{M^y F}{G} =
\int_{{X}_{\mathbb R}}
f(x)y(x) \overline{ g(x)}
\frac{d \eta(x)}{\omega(x)}
\]
and from the assumption that $y$ is a real meromorphic function.
\end{pf}

Equation \ref{MyLPhi} immediately produces the $\mathcal{L}(\varphi)$ counterpart of \cite[Theorem 4.6]{av3}.

\begin{corollary}
Let $y_1$ and $y_2$ be two meromorphic functions of degree $n_1$ and $n_2$, respectively.
Furthermore, we assume that the poles of $y_1$ and $y_2$ lie outside the support of $d \eta$ on $X _{\mathbb R}$.
Then, $M^{y_1}$ and $M^{y_2}$ commute on $\mathcal{L}(\varphi)$, that is, for every
$F(z) \in  \mathcal{L}(\varphi)$, $M^{y_1} M^{y_2} F = M^{y_2} M^{y_1} F$ holds.
\end{corollary}

Using the observation that $R^y_{\alpha}$ is just the operator $M^{\frac{1}{y(u)-\alpha}}$, we present the counterpart of \eqref{tatche},
that is, an integral representation of the resolvent operator at $\alpha$.

\begin{corollary}
$\mathcal{L}(\varphi)$ is invariant under the resolvent operator $R^y_\alpha$ where $\alpha$ is a non-real complex number.
Moreover, the resolvent operator has the integral representation:
\begin{equation}
\label{RyInLPhi}
(R^y_\alpha F)(u)
=
\int_{X _ {\mathbb R}}K_{\zeta}(u,x)\frac{f(u)}{(y(u)-\alpha)} \frac{d \eta(x)}{\omega(x)}
,
\end{equation}
where the poles of $y$ do not belong to the support of $d \eta$.
\end{corollary}

As another immediate corollary,
we mention that any pair of resolvent operator commutes.

\begin{corollary}
Let $y_1$ and $y_2$ be two meromorphic functions of degree $n_1$ and $n_2$, respectively.
Furthermore, assume the poles of $y_1$ and $y_2$ lie outside the support of $d \eta$ on $X _{\mathbb R}$ and
let $\alpha$ and $\beta$ be two elements in $\mathbb C \setminus \mathbb R$.
Then the resolvent operators $R^{y_1}_{\alpha}$ and $R^{y_2}_{\beta}$ commute,
i.e. for any $F(z) \in  \mathcal{L}(\varphi)$ the following holds $R^{y_1}_{\alpha} R^{y_2}_{\beta} F= R^{y_2}_{\beta} R^{y_1}_{\alpha} F$.
\end{corollary}
The counterpart of Theorem \ref{finiteDimentionalLphi} is given below.

\begin{Tm}
Let $\varphi$ be analytic in $X \backslash X_{\mathbb R}$ and with positive real part.
Then the following are equivalent:
\begin{enumerate}
\item \label{tfre1} The reproducing kernel space $\mathcal{L}(\varphi)$ is finite dimensional.
\item \label{tfre2} $\varphi$ is meromorphic on ${X}$.
\item \label{tfre3} ${\bf L}^2(d \eta)$ is finite dimensional.
\end{enumerate}
\end{Tm}
\begin{pf}
Since $\mathcal{L}(\varphi)$ is isomorphic to ${\bf L}^2(d \eta)$,
$\mathcal{L}(\varphi)$ is finite dimensional 
if and only if ${\bf L}^2(d \eta)$ is finite dimensional.
\smallskip

If ${\bf L}^2(d \eta)$ is finite dimensional then $d \eta$ has a finite number of atoms and hence $\varphi$ is meromorphic.
On the other hand assume that $\varphi$ is meromorphic on $X$.
As in the classical case, the measure in the integral representation of $\varphi$ is obtained as the weak star limit of $\varphi(p) + \overline{\varphi(p)}$,
and the poles of $\varphi$ correspond to the atoms of $d \eta$.
\end{pf}

\begin{Tm}
\label{thm432}
Let $f,g \in \mathcal{L}(\varphi)$ and $\alpha, \beta \in \mathbb C$ with non-zero imaginary part.
Then the following identity holds,
\begin{equation}
\label{strucIdPhi}
\innerProductReg{R^y_\alpha f}{ g} -
\innerProductReg{f}{ R^y_\beta g} -
(\alpha - \overline{\beta}) \innerProductReg{R^y_\alpha f}{ R^y_\beta g} = 0.
\end{equation}
\end{Tm}
\begin{pf}
Using \eqref{RyInLPhi}, the left hand side of \eqref{strucIdPhi} can be written as
\begin{align*}
\innerProductReg{R^y_\alpha f}{ g} & -
\innerProductReg{f}{ R^y_\beta g} -
(\alpha - \overline{\beta}) \innerProductReg{R^y_\alpha f}{ R^y_\beta g}
=
\int_{X _ {\mathbb R}}
f(u) g(u) K_{\zeta}(p,u)
\times
\\
&
\left(
\frac{1}{(y(u)-\alpha)}
-
\frac{1}{(y(u)-\overline{\beta})}
-
\frac{\alpha - \overline{\beta}}{(y(u)-\alpha)(y(u)-\overline{\beta})}
\right)
\frac{ d \eta(u)}{\omega(u)}.
\end{align*}

One may note that
\[
\frac{1}{(y(u)-\alpha)}
-
\frac{1}{(y(u)-\overline{\beta})}
-
\frac{\alpha - \overline{\beta}}{(y(u)-\alpha)(y(u)-\overline{\beta})}
\]
is identically zero, hence the result follows.

\end{pf}

In fact Theorem \ref{thm432} is an if and only if relation, and we refer the reader to \cite{AVP3}
for the related de Branges structure theorems.


\section{The \texorpdfstring{$\mathcal{L}(\varphi)$}{ $\mathcal{L}(\varphi)$ } spaces in the single-valued case}
\label{secPhiSingleVal}

Whenever an additive function $\varphi$ is single-valued,
the formula
\begin{equation*}
s(p)= \frac{1-\varphi(p)}{1+\varphi(p)}
\end{equation*}
makes sense and defines a single-valued
function $s(p)$.
Then, the reproducing kernel associated with $s(p)$, denoted by ${\mathcal H}(s)$, 
is of the form
\[
i(1-s(p)s(q)^*)K_{\zeta}( p,{\tauBa{q}}).
\]
These spaces were studied in \cite{av3} in the finite dimensional setting and
in \cite{AVP1} in the infinite dimensional case.
\smallskip

We note that the multiplication operator $ u \mapsto \frac{(1+\varphi(p))}{{\sqrt{2}}} u $ maps, 
as in the zero genus case, ${\mathcal H}(s)$ onto $\mathcal{L}(\varphi)$ unitarily.
Hence, we may pair any $u\in{\mathcal H}(s)$ to a function
$f\in{\bf L}^2(d \eta)$ through the corresponding $\mathcal{L}(\varphi)$ space, such that
\begin{equation*}
\frac{1}{\sqrt{2}}(1+\varphi(p))u(p)
=
\int_{{X}_{\mathbb R}} K_{\zeta}( p,x)f(x)\frac{d \eta(x)}{\omega(x)}.
\end{equation*}
We denote the mapping from ${\mathcal H}(s)$ onto ${\bf L}^2(d \eta)$ by $$\Lambda: u(p) \longrightarrow f(x).$$
We now turn to express the operator $M^y$ using the operator of multiplication by $y$ in
${\bf L}^2(d \eta)$.

\begin{Tm}
Let $\varphi$ be a single-valued function with positive real part on a dividing-type compact Riemann surface $X$
and let $y$ be a meromorphic function of degree $n$ on $X$.
Then, any $f\in {\bf L}^2(d \eta)$ satisfies
\begin{align}
\nonumber
(\Lambda M^y \Lambda^*)f(x) = &
y(x)f(x)+
\\
\label{jardin-des-plantes}
&
i
\sum_{j=1}^{n}
\frac{c_j}{1+\varphi(p^{(j)})}K_{\zeta}( x, p^{(j)})
\left(\int_{{X}_{\mathbb R}} K_{\zeta}( p^{(j)}, p)f(p)
\frac{d \eta(p)}{\omega(p)}\right),
\end{align}
where $p^{(1)},...,p^{(n)}$ are the $n$ distinct poles of $y$.
\end{Tm}
\begin{pf}
Let $u\in{\mathcal H}(s)$ and $f\in{\bf L}^2(d \eta)$ such that, $\Lambda u = f$, that is,
they satisfy the relation
\begin{equation}
\label{phiInt}
\frac{1+\varphi(p)}{{\sqrt{2}}}u(p)=\int_{{ X}_{\mathbb
R}}K_{\zeta}(x,p)f(x)\frac{d \eta(x)}{\omega(x)}.
\end{equation}
Then, multiplying both sides of \eqref{phiInt} by $y(p)$, we obtain
\begin{align}
\nonumber
y(p)
\frac{1+\varphi(p)}{\sqrt{2}}u(p)
= &
\int_{{ X}_{\mathbb R}}y(p)K_{\zeta}( p, x)f(x)\frac{ d \eta(x)}{\omega(x)}
\\
\nonumber
=&
\int_{{ X}_{\mathbb R}}y(x)K_{\zeta}( p, x)f(x)\frac{d \eta(x)}{\omega(x)}+
\\
&
\int_{{ X}_{\mathbb R}}(y(p)-y(x))K_{\zeta}( p, x)f(x) \frac{d \eta(x)}{\omega(x)}.
\label{phiInt1}
\end{align}
Then, using the collection-type formula (see \cite[Proposition 3.1, Eq. 3.5]{av2}), we have
\begin{equation*}
(y(p)-y(q))K_{\zeta}( p, q)
=
-
\sum_{j=1}^{n}
\frac{c_j}{dt_j(p^{(j)})}K_{\zeta}( p, p^{(j)})K_{\zeta}( p^{(j)},q).
\end{equation*}
Then \eqref{phiInt1} becomes:
\begin{align}
\nonumber
\int_{{X}_{\mathbb R}}y(x)K_{\zeta}( p, x)f(x)
\frac{d \eta(x)}{\omega(x)}
=
&
\frac{(1+\varphi(p))}{\sqrt{2}}
y(p)u(p) -
\\
&
\sum_{j=1}^{n}
\frac{c_j K_{\zeta}( p, p^{(j)})}{dt_j(p^{(j)})}
\int_{{X}_{\mathbb R}} K_{\zeta}( p^{(j)}, x)f(x)\frac{d \eta(x)}{\omega(x)}.
\label{phiInt2}
\end{align}
On the other hand, using the equality
\begin{equation}
(1+\varphi(p))\displaystyle{\frac{1+s(p)}{2}}
=
1.
\label{phiInt3}
\end{equation}
Equation \ref{phiInt} becomes
\begin{equation}
\label{phiInt4}
\int_{{X}_{\mathbb R}}
K_{\zeta}(x,p^{(j)})f(x)\frac{d \eta(x)}{\omega(x)}=
\displaystyle{\frac{\sqrt{2}}{1+s(p^{(j)})}}u(p^{(j)}).
\end{equation}
Now, substituting \eqref{phiInt3} and \eqref{phiInt4} in \eqref{phiInt2},
we obtain the following calculation:
\begin{align}
\frac{\sqrt{2}}{1+\varphi(p)}
&
\int_{{X}_{\mathbb R}}
K_{\zeta} ( p, x) y(x)f(x) \frac{d \eta(x)}{\omega(x)} =
\nonumber
\\
= &
y(p)u(p)
-
\frac{1+s(p)}{\sqrt{2}}
\sum_{j=1}^{n}
c_j K_{\zeta}( p, p^{(j)})\frac{\sqrt{2}u(p^{(j)})}{1+s(p^{(j)})}
\nonumber
\\
=&
y(p)u(p)
-
\sum_{j=1}^{n}
c_j\frac{1+s(p)}{1+s(p^{(j)})}
K_{\zeta}( p , p^{(j)})u(p^{(j)})
\nonumber
\\
\nonumber
=&
y(p)u(p)
-
\sum_{j=1}^{n}
c_j
K_{\zeta}( p , p^{(j)})u(p^{(j)})
\left(1 + \frac{s(p)-s(p^{(j)})}{1+s(p^{(j)})} \right)
\\
=
&
(M^y u)(p)
+
\sum_{j=1}^{n}
c_j
\frac{s(p)-s(p^{(j)})}{\sqrt{2}}
K_{\zeta}( p , p^{(j)})
u(p^{(j)})
.
\label{phiInt5}
\end{align}
On the other hand, using \eqref{phiInt3}, we have
\begin{align}
\nonumber
\varphi(p)-\varphi(p^{(j)})
=&
\frac{2(s(p^{(j)})-s(p))}{(1+s(p))(1+s(p^{(j)}))}
\\
=&
(1+\varphi(p))(1+\varphi(p^{(j)}))
\frac{s(p^{(j)})-s(p)}{2}.
\label{phiInt6}
\end{align}
Thus, we take \eqref{phiInt6} and multiply it on the right by the Cauchy kernel $K_{\zeta}( p, p^{(j)})$ in both sides
and use \eqref{jfk-le-26-octobre-2000_SV} to conclude
\begin{align*}
(1+\varphi(p))
\frac{s(p)-s(p^{(j)})}{2}
&
K_{\zeta}( p, p^{(j)})
= 
\frac{\varphi(p^{(j)})-\varphi(p)}{1+\varphi(p^{(j)})}
K_{\zeta}( p, p^{(j)})
\\
=&
\frac{i}{1+\varphi(p^{(j)})}
\int_{{X}_{\mathbb R}}
K_{\zeta}( p, x)K_{\zeta}( x,p^{(j)})\frac{d \eta(x)}{\omega(x)}
.
\end{align*}
Thus, \eqref{phiInt5} becomes
\begin{align*}
\frac{1+\varphi(p)}{\sqrt{2}} & (M^yu)(p)
= 
\int_{{X}_{\mathbb R}}
K_{\zeta}( p, x)y(x)f(x)
\frac{d \eta(x)}{\omega(x)}
+
i \sum_{j=1}^{n}
\frac{c_j}{1+\varphi(p^{(j)})} \times
\nonumber
\\
&
\times\left(
\int_{{X}_{\mathbb R}}K_{\zeta}( p, x)K_{\zeta}( x,p^{(j)})\frac{d \eta(x)}{\omega(x)}
\right)
\left(
\int_{{X}_{\mathbb R}} K_{\zeta}( p^{(j)}, s)f(s)\frac{d \eta(s)}{\omega(s)}
\right),
\end{align*}
and by setting
\begin{align*}
{\bf \widehat{f}}(q)
= 
y(q)f(q)
+
i
\sum_{j=1}^{n}
\frac{c_jK_{\zeta}( q, p^{(j)})}{1+\varphi(p^{(j)})}
\left(
\int_{{X}_{\mathbb R}} K_{\zeta}( p^{(j)}, x)f(x) \frac{d \eta(x)}{\omega(x)}
\right),
\end{align*}
the identity in \eqref{phiInt6} becomes
\[
\frac{1+\varphi(p)}{\sqrt{2}}(M^yu)(p)
=
\int_{{X}_{\mathbb R}} K_{\zeta}( p, x){\bf \widehat{f}}(x)\frac{d \eta(x)}{\omega(x)}
.
\]
\end{pf}
\begin{Cy}
Let $\varphi$ be a single-valued function with positive real part on a dividing-type compact Riemann surface $X$.
Furthermore, let us assume that $y(p)$ is a real meromorphic function of degree $n$ such that $s(p^{(j)})=1$ for all $1\leq j \leq n$.
Then the following identity holds:
\begin{equation*}
(\Lambda (\real M^y)\Lambda^*)f(p)=y(p)f(p).
\end{equation*}
\end{Cy}
We note that, in fact, one may assume that $s(p^{(j)})$ for all $1 \leq j \leq n$ equal to a
common constant of modulus one.
\smallskip

Furthermore, for an arbitrary $f\in{\bf L}^2(d \eta)$ we set
\[
\Phi_{y}(f)
\defEq
\col_{1\leq j \leq n}~\left(
\frac{1}{1+\varphi(p^{(j)})}\int_{{X}_{\mathbb R}}
K_{\zeta}( p^{(j)}, x)f(x)\frac{d \eta(x)}{\omega(x)}\right),
\]
and then, for an element $d= (d_1,...,d_n)^t \in {\mathbb C}^{n}$, the adjoint operator
$\Phi_{y}^*:\C^n \rightarrow \mathcal{L}(\varphi)$ is given explicitly by
\[
\Phi_{y}^*d
=
\sum_{j=1}^{n} \frac{1}{1+\overline{\varphi(p^{(j)})}}d_j K_{\zeta}( p, \tauBa{p^{(j)}}).
\]
Then, if we further use the notation
$$\sigma_y \defEq {\rm diag}~ c_j \, (1+\overline{\varphi(p^{(j)})}),$$
the formula in \eqref{jardin-des-plantes} may be rewritten and simplified as follows
\begin{equation*}
(\Lambda M^y \Lambda^*)f(p)
=
y(p)f(p)+\frac{i}{2}\Phi_{y}^*\sigma_y\Phi_{y} f,
\end{equation*}
while the real part of the operator $M^y$ has the form
\begin{equation*}
(\Lambda (\real M^y)\Lambda^*)f(p)=y(p)f(p)+
\Phi_{y}^*{\rm diag}~({\rm Im}~\varphi(p^{(j)}))\Phi_{y} f.
\end{equation*}

\begin{landscape}

\section{Summary}
\label{chSumm43}
\setcounter{equation}{0}
The table below summarizes the comparison between the $\mathcal{L}(\varphi)$ spaces 
in the Riemann sphere case and in a compact real Riemann surfaces setting.
\begin{center}
\begin{tabular}{|m{5cm}||M{6.0cm}|M{8.3cm}|}
\hline
&
{\bf The $g=0$ setting}
&
{\bf The $g>0$ setting}
\\
\hline\hline
&
$z-w$
&
$E(p,q)$
\\
\hline
The Cauchy kernel
&
$\frac{1}{-i(z-\overline{w})}$
&
$K_{\zeta}(p,\tauBa{q})
\defEq
\frac
{\vartheta[\zeta](p-{\tauBa{q}})} {i\vartheta[\zeta](0)E(p,{\tauBa{q}})}
$
\\
\hline
The Hardy space $H^2$
&
The reproducing kernel Hilbert space with kernel $\frac{1}{-i(z-\overline{w})}$
&
The reproducing kernel Hilbert space the with kernel
$\frac{1}{-i}K_{\zeta}(p, \tauBa{q})$
where $\zeta\in T_0$.
\\
\hline
The kernel of $\mathcal{L}(\varphi)$
&
$\frac{\varphi(z)+\varphi(w)^*}{-i(z-\overline{w})} = \innerProductTri{\frac{1}{t-z}}{\frac{1}{t-w}}{{\bf L}^2(d \eta)},$
&
$ \left({\varphi}(p)+{\varphi}(q)^*\right)
\frac{\vartheta[\zeta](p-{\tauBa{q}})}{\vartheta[\zeta](0)E({\tauBa{q}},p)}$
\\
\hline
Reproducing kernel in ${\bf L}^2( d \mu) $
&
$\int_{\mathbb R}\frac{d \eta(t)}{(t-z)(t-\overline{w})}$
&
\small
$\begin{aligned} 
\innerProductTri {K_{\zeta}(\tauBa{q},x)}{ K_{\zeta}( \tauBa{p}, x)}{{\bf L}^2\left(\frac{d \eta}{ \omega}\right)} 
\end{aligned}$
\\
\hline
The elements of $\mathcal{L}(\varphi)$
&
$F(z)=\int_{\mathbb R}\frac{f(t)d \eta(t)}{t-z}$
&
\small
$
\innerProductReg {f}{ K_{\zeta}( \tauBa{p}, x)} 
=
\int_{{X}_{\mathbb R}}f(x)K_{\zeta}(x, p)f(x)\frac{d \eta(x)}{\omega(x)}
$
\\
\hline
{The Herglotz integral representation formula}
&
\small
$\begin{aligned} 
\varphi(z)= & iA-iBz + 
\\ &
i\int_{{\mathbb R}}\left(\frac{1}{t-z}-\frac{t}{t^2+1} \right)d \eta(t)
\end{aligned}$
&
\small
$\setlength{\jot}{0pt}\begin{aligned}
\varphi(z) = & 
\frac{\pi}{2} \int_{X_{\mathbb R}} \frac{[\omega_1(x)  \cdots  \omega_g(x)]}{\omega(x)} n(\widetilde{x}) \, d \eta(x) +
\\ &
\pi i\int_{X_{\mathbb R}} \frac{[\omega_1(x)\,\cdots \,\omega_g(x)]}{\omega(x)}(Yp)d \eta(x) -
\\ &
\frac{i}{2}\int_{X_{\mathbb R}}\frac{ \frac{\partial}{\partial x} \ln E(p,\widetilde{x})}{\omega(x)}d \eta(x) + iM
\end{aligned}$
\\
\hline
Integral representation of the model operator $M^y$
&
$(M F)(z)=\int_{\mathbb R}{\frac{t f(t)d \eta(t)}{t-z}}$
&
$
(M^y F)(z)=
\int_{{X}_{\mathbb R}}
K_{\zeta}( u, x)
f(x)y(x)\frac{d \eta(x)}{\omega(x)}$
\\
\hline
Integral representation of the resolvent operator $R_\alpha^y$
&
$(R_\alpha F)(z)=\int_{\mathbb R}\frac{f(t)d \eta(t)}{(t-z)(t-\alpha)}$
&
\small
$\begin{aligned} 
(R_\alpha F)(p) & =
\int_{ X_{\mathbb R}}
K_{\zeta}( p,u) \frac{f(u)}{y(u) - \alpha}
\frac{ d \mu(u)}{\omega(u)}
\end{aligned}$
\\
\hline
\end{tabular}
\end{center}
\vspace{0.5cm}
\end{landscape}
\normalsize


\def\cfgrv#1{\ifmmode\setbox7\hbox{$\accent"5E#1$}\else
  \setbox7\hbox{\accent"5E#1}\penalty 10000\relax\fi\raise 1\ht7
  \hbox{\lower1.05ex\hbox to 1\wd7{\hss\accent"12\hss}}\penalty 10000
  \hskip-1\wd7\penalty 10000\box7} \def\cprime{$'$} \def\cprime{$'$}
  \def\cprime{$'$} \def\lfhook#1{\setbox0=\hbox{#1}{\ooalign{\hidewidth
  \lower1.5ex\hbox{'}\hidewidth\crcr\unhbox0}}} \def\cprime{$'$}
  \def\cprime{$'$} \def\cprime{$'$} \def\cprime{$'$} \def\cprime{$'$}
  \def\cprime{$'$}


\end{document}